\documentclass[reqno,a4wide,english]{amsart}

\usepackage{amsmath}

\usepackage[all]{xy}

\theoremstyle{plain}
\newtheorem{theorem}{Theorem}[section]

\theoremstyle{definition}

\newtheorem{emptythm}[theorem]{}

\newcommand{\too}{\longrightarrow}

\newcommand{\smb}{{\scriptscriptstyle \bullet}}

\newcommand{\ol}[1]{\overline{#1}}
\newcommand{\ul}[1]{\underline{#1}}

\DeclareMathOperator{\Hom}{Hom}

\DeclareMathOperator{\Ker}{Ker}
\DeclareMathOperator{\Coker}{Coker}
\DeclareMathOperator{\id}{id}
\DeclareMathOperator{\rank}{rank}

\DeclareMathOperator{\khar}{char}
\DeclareMathOperator{\Spec}{Spec}

\DeclareMathOperator{\Br}{\mathrm{Br}}
\DeclareMathOperator{\Pic}{\mathrm{Pic}}
\DeclareMathOperator{\HM}{\mathrm{H}}
\DeclareMathOperator{\et}{\mathrm{et}}
\DeclareMathOperator{\sico}{\mathrm{sc}}
\DeclareMathOperator{\Gal}{\mathrm{Gal}}
\DeclareMathOperator{\aff}{\mathrm{aff}}
\DeclareMathOperator{\CH}{\mathrm{CH}}
\DeclareMathOperator{\Gl}{\mathrm{GL}}
\DeclareMathOperator{\Sl}{\mathrm{SL}}
\DeclareMathOperator{\Rad}{\mathrm{Rad}}
\DeclareMathOperator{\rk}{\mathrm{rk}}

\DeclareMathOperator{\res}{\mathrm{res}}
\DeclareMathOperator{\cor}{\mathrm{cor}}

\newcommand{\MK}[1]{\mathrm{K}^{M}_{#1}}  
\newcommand{\mk}[1]{\mathrm{k}^{M}_{#1}}  

\DeclareMathOperator{\CK}{\mathrm{C}}     

\newcommand{\CM}[1]{\mathrm{M}_{#1}}     

\newcommand{\bBr}{\overline{\Br}}

\newcommand{\Brtor}[1]{{ }_{#1}\hspace{-0.5mm}\Br}
\newcommand{\HMtor}[1]{{ }_{#1}\hspace{-0.5mm}\HM}

\newcommand{\roun}{\mu}                                
\newcommand{\ru}{\mu\hspace{-1.5mm}\mu}                
\newcommand{\Gm}{\mathbb{G}_{{\rm m}}}                  

\newcommand{\Z}{\mathbb{Z}}
\newcommand{\N}{\mathbb{N}}
\renewcommand{\P}{\mathbb{P}}
\newcommand{\A}{\mathbb{A}}

\newcommand{\maxid}{\mathfrak{m}}

\newcommand{\cf}{\textsl{cf.}\ }
\newcommand{\eg}{\textsl{e.g.}\ }
\newcommand{\ie}{\textsl{i.e.}\ }

\begin{document}

\title[On the Brauer group of a product]
{On the Brauer group of the product of a torus and a semisimple algebraic group}

\author{Stefan Gille}
\email{gille@ualberta.ca}
\address{Department of Mathematical and Statistical Sciences,
University of Alberta, Edmonton T6G 2G1, Canada}

\author{Nikita Semenov}
\email{semenov@uni-mainz.de}
\address{Johannes Gutenberg-Universit\"at Mainz, Institut f\"ur
Mathematik, Staudingerweg 9, 55099 Mainz, Germany}

\thanks{The first named author was partially supported
by an NSERC research grant. The authors acknowledge a partial support of SFB/Transregio 45 
``Periods, moduli spaces and arithmetic of algebraic varieties''.}


\date{January 24, 2013}

\begin{abstract}
Let~$T$ be a torus (not assumed to be split) over a field~$F$,
and denote by~$\HMtor{n}^{2}_{\et}(X,\Gm)$ the subgroup of elements
of exponent dividing~$n$ in the cohomological Brauer group of a
scheme~$X$ over the field~$F$. We provide conditions on $X$ and~$n$
for which the pull-back homomorphism
$\HMtor{n}^{2}_{\et}(T,\Gm)\too\HMtor{n}^{2}_{\et}(X\times_{F}T,\Gm)$
is an isomorphism. We apply this to compute the Brauer group of
some reductive groups and of non singular affine quadrics.

Apart from this, we investigate the $p$-torsion of the Azumaya
algebra defined Brauer group of a regular affine scheme over
a field~$F$ of characteristic~$p>0$.
\end{abstract}

\maketitle

\section{Introduction}
\label{IntroSect}\bigbreak

\noindent
Let~$G$ be a simply connected group over an algebraically
closed field~$F$ of characteristic zero. Magid~\cite{Ma78}
has shown that the pull-back homomorphism
$$
p_{T}^{\ast}\colon\;\Br (T)\,\too\,\Br (G\times_{F}T)
$$
along the projection $p_{T}\colon G\times_{F}T\too T$
is an isomorphism for all $F$-tori~$T$. Here and in
the following we denote by~$\Br (X)$ the Brauer group
of equivalence classes of Azumaya algebras over a
scheme~$X$, and by~$\Brtor{n}(X)$ the subgroup of elements
whose exponent divides the integer~$n\geq 2$.

\smallbreak

Magid's proof uses the Lefschetz principle to reduce to
the case that~$F$ is the field of complex numbers, and
then Grothendieck's deep comparison result for \'etale
and singular cohomology with finite coefficients, which
finally reduces everything to computations of singular
cohomology groups.

\medbreak

In the present article we give a purely algebraic proof
of Magid's result and generalize it to groups
over arbitrary fields (see Corollary~\ref{SemisimpleGrThm} below).
We also generalize the main theorem of~\cite{Gi09}
which asserts that the natural homomorphism
$\Brtor{n}(F)\too\Brtor{n}(G)$ is an isomorphism for
any semisimple algebraic group~$G$ over a field~$F$
if~$\khar F$ does not divide~$n$ and~$n$ is
coprime to the order of the fundamental group of~$G$.

\smallbreak

Actually, we prove a more general result,
see Theorem~\ref{mainThm}. Let~$\HM^{2}_{\et}(X,\Gm)$
be the cohomological Brauer group of the scheme~$X$,
and denote by~$\HMtor{n}^{2}_{\et}(X,\Gm)$ the subgroup
of elements whose order divides~$n$. We give conditions
under which for a regular and integral
scheme~$X$ of finite type over a field~$F$
the natural homomorphisms
$$
\Brtor{n}(F)\,\too\,\HMtor{n}^{2}_{\et}(X,\Gm)\qquad
\mbox{and}\qquad
\Brtor{n}(T)\too\HMtor{n}^{2}_{\et}(T\times X,\Gm)\, ,
$$
where $T$ is a form of a torus, are isomorphisms. (Recall
that by a result of Hoobler~\cite{Ho80} the natural homomorphism
$\Br (X)\too\HM^{2}_{\et}(X,\Gm)$ defined in the work~\cite{Gr68}
by Grothendieck is an isomorphism for all smooth schemes~$X$
over a field.)

\smallbreak

As an application we get from Theorem~\ref{mainThm}
some computations of Brauer groups
of reductive groups and affine quadrics. For instance if~$\khar F=0$
we show in Corollary~\ref{cor312} that if~$H$ is a reductive group,
such that the derived group~$G=[H,H]$ is adjoint,
then $\Brtor{n}(H)\simeq\Brtor{n}(\Rad (H))$
for all integers~$n\geq 2$ which are coprime
to the order of the fundamental group of~$G$.
(Here $\Rad(H)$ stands for the radical of $H$
and~$[H,H]$ denotes the commutator subgroup of~$H$.)

\smallbreak

Theorem~\ref{mainThm} also implies that if~$X$ is a
non singular affine quadric of Krull dimension~$\geq 3$
over a field~$F$ of characteristic~$\not= 2$
then the pull-back of the structure morphism
$$
\Brtor{n}(F)\,\too\,\Brtor{n}(X)
$$
is an isomorphism for all integers~$n$ which are not
divisible by~$\khar F$. In particular we have an isomorphism
$\Br (F)\xrightarrow{\simeq}\Br (X)$ for such affine
quadrics over a field~$F$ of characteristic zero.

\medbreak

In the last section we investigate the $p$-torsion of the Brauer group,
where $p$ is the characteristic of the base field. We show that if~$X$ is
the spectrum of a regular geometrically integral algebra of
finite type and of Krull dimension~$\geq 2$ over a field~$F$
of characteristic~$p>0$ then
the cokernel of the pull-back $\Br (F)\too\Br (X)$ contains a
non empty $p$-torsion and $p$-divisible (hence infinite)
abelian group. In particular if~$G$ is a linear
algebraic group of Krull dimension~$\geq 2$ over a field~$F$
of characteristic~$p>0$ then the pull-back $\Br (F)\too\Br (G)$
is never surjective.

\smallbreak

This gives an explanation that our computational results
for the $n$-torsion subgroup of the Brauer group of
certain affine schemes over a field~$F$ are only true if~$\khar F$
does not divide~$n$.

\smallbreak

Note that this is an ``affine'' result. For a geometrically
integral smooth affine variety~$X$ of dimension~$\geq 2$
over a field~$F$ of positive characteristic the natural
morphism $\Br (F)\too\Br (X)$ is never surjective, but
there are a examples of geometrically integral
smooth and projective varieties~$Y$ over~$F$, such that
$\Br (F)\too\Br (Y)$ is surjective (and hence an isomorphism
if~$Y$ has an $F$-rational point). For instance, Merkurjev
and Tignol~\cite[Thm.\ B]{MeTi95} have shown that for a
twisted flag variety~$Y$ for a semisimple algebraic
group~$G$ over an arbitrary field~$F$ the homomorphism
$\Br (F)\too\Br (Y)$ is surjective.

\bigbreak\bigbreak

\noindent
{\bfseries Acknowledgment.}
The first author would like to thank Vladimir Chernousov,
Jochen Kuttler, and Alexander Vishik for discussions and comments,
and the Institut f\"ur Mathematik der Johannes Gutenberg-Universit\"at Mainz for
hospitality during a two weeks stay in July 2012. Both authors
would like to thank Victor Petrov for discussions and comments.

\bigbreak

\goodbreak
\section{Preliminaries on the Brauer group}
\label{PreBrauerGroupSect}\bigbreak

\begin{emptythm}
\label{NotationSubSect}
{\it Notation.}
Let~$F$ be a field with separable closure~$F_{s}$.
We denote by~$\Gamma_{F}=\Gal (F_{s}/F)$ the absolute
Galois group of~$F$. If~$X$ is an $F$-scheme we set
$X_{s}:=F_{s}\times_{F}X$.

\smallbreak

By $\HM^{i}(F,\, -\, )$ we denote the Galois cohomology
of the field~$F$.

\smallbreak

We denote by~$R^{\times}$ the multiplicative
group of units of a commutative ring~$R$.

\smallbreak

Let~$X$ be an $F$-scheme. Then we denote by~$F[X]$
its ring of regular functions, and in case~$X$ is integral
by~$F(X)$ its function field.

\bigbreak

For a prime number~$p$ we denote by~$\Z (p^{\infty})$
the infinite abelian group $\Z [\frac{1}{p}]/\Z$. This
is a $p$-divisible and $p$-torsion group and any
$p$-divisible and $p$-torsion group is a direct
sum of copies of~$\Z (p^{\infty})$, see \eg~\cite{IAb}.
%
%
\end{emptythm}

\begin{emptythm}
\label{BrauerGrSubSect}
{\it Brauer groups.}
Let~$R$ be a commutative ring (with~$1$) and
$X=\Spec R$. The Brauer group~$\Br (R)$ of the
commutative ring~$R$ (also denoted $\Br (X)$) has been
introduced by Auslander and Goldman~\cite{AuGo60}.
It classifies Azumaya algebras over~$R$ and generalizes
the classical Brauer group of a field.

\smallbreak

If~$A$ is an $R$-Azumaya algebra we denote by~$[A]$
the class of~$A$ in~$\Br (R)$.

\medbreak

There is also a cohomological Brauer group of~$R$, the second
\'etale cohomology group~$\HM^{2}_{\et}(X,\Gm)$. By
Grothendieck~\cite{Gr68} we have a natural
homomorphism
$$
j_{X}\colon\;\Br (X)\,\too\,\HM^{2}_{\et}(X,\Gm)\, .
$$
If~$X=\Spec R$ is smooth over a field~$F$, then~$j_{X}$
is an isomorphism as shown by Hoobler~\cite{Ho80}. In
this case we identify~$\Br (X)$ with~$\HM^{2}_{\et}(X,\Gm)$.
(There is a more general identification result due
to Gabber but we do not need this.)

\medbreak

We denote the subgroup of elements of order dividing~$n$ in~$\Br (X)$
by~$\Brtor{n}(X)$, and similar we denote by~$\HMtor{n}^{2}_{\et}(X,\Gm)$
the subgroup of elements of order~$n$ in the cohomological Brauer
group of~$X$.

\smallbreak

If~$\khar F$ does not divide~$n$ we have an exact sequence
of \'etale sheaves
\begin{equation}
\label{cdot-nEq}
\xymatrix{
{\ru_{n}}\;\, \ar@{>->}[r] & {\Gm} \ar@{->>}[r]^-{\cdot n} & {\Gm}\, ,
}
\end{equation}
\ie~$\ru_{n}$ is (by definition) the kernel of $\Gm\xrightarrow{\cdot n}\Gm$.
It follows from the associated long exact cohomology sequence
that the homomorphism
$$
\HM^{2}_{\et}(X,\ru_{n})\,\too\,\HMtor{n}^{2}_{\et}(X,\Gm)
$$
is an isomorphism if~$\Pic X$ is $n$-divisible (\eg if
$\Pic X$ is trivial).
\end{emptythm}

\begin{emptythm}
\label{HSSpSeqSubSect}
{\it The Hochschild-Serre spectral sequence.}
This spectral sequence is a useful tool to compute
Brauer groups. Let $X'\too X$ be a Galois
cover with (finite) group~$\Pi$. Then,
see~\cite[Chap.\ III, Thm.\ 2.20]{EC}, there is
a convergent spectral sequence
$$
E_{2}^{p,q}\, :=\;\HM^{p}(\Pi,\HM^{q}_{\et}(X',\Gm))\;\;\;
\Longrightarrow\;\HM^{p+q}_{\et}(X,\Gm)\, ,
$$
where $\HM^{i}(\Pi,\, -\, )$ denotes the group cohomology of~$\Pi$.

\smallbreak

If~$X$ is a finite type $F$-scheme we get from this spectral
sequence for finite Galois covers by a limit argument,
see~\cite[Rem.\ 2.21]{EC}, the following spectral sequence
\begin{equation}
\label{HSSpSeqEq}
E_{2}^{p,q}\, :=\;\HM^{p}(F,\HM^{q}_{\et}(X_{s},\Gm))\;\;\;
\Longrightarrow\;\HM^{p+q}_{\et}(X,\Gm)\, ,
\end{equation}
where (as in~\ref{NotationSubSect}) we have set
$X_{s}=F_{s}\times_{F}X$ and~$\HM^{i}(F,\, -\,)$ denotes
Galois cohomology of the field~$F$.

\smallbreak

For later use we recall that if~$\Pic X_{s}=0$ and~$X$ is a
smooth and affine $F$-scheme we get from this spectral sequence
an exact sequence
$$
0\too\HM^{2}(F,F_{s}[X_{s}]^{\times})\too\HM^{2}_{\et}(X,\Gm)\too
\HM^{2}_{\et}(X_{s},\Gm)^{\Gamma_{F}}\too\HM^{3}(F,F_{s}[X_{s}]^{\times})\, ,
$$
where~$\Gamma_{F}$ is the absolute Galois group of the field~$F$.
\end{emptythm}

\begin{emptythm}
\label{AzExplSubSect}
{\it A class of examples of ``cyclic'' Azumaya algebras.}
Let~$R$ be an integral domain
and~$x,y\in R^{\times}$. Let further $n\geq 2$
be an integer and assume that~$R$ contains~$\frac{1}{n}$ and
a primitive $n$-th root of unity~$\xi$.

\smallbreak

We define an associative $R$-algebra~$A=A_{\xi}^{R}(x,y)$
as follows: The~$R$-algebra~$A$ is generated by elements~$\alpha$
and~$\beta$ which are subject to the following relations:
$$
\alpha^{n}=x\, ,\;\,\beta^{n}=y\, ,\;\mbox{and}\;\;
x\cdot y=\xi (y\cdot x)\, .
$$
If~$\maxid$ is a maximal ideal of the ring~$R$, then by
Milnor~\cite[Chap. 15]{IalgKT} we know that
$R/\maxid\otimes_{R} A_{\xi}^{R}(x,y)$ is a central
simple $R/\maxid$-algebra and so~$A_{\xi}^{R}(x,y)$
is an $R$-Azumaya algebra.

\smallbreak

Assume now that~$R$ is a regular integral domain. Let~$K$ be
its quotient field. Then the natural morphism $\Br (R)\too\Br (K)$
is injective as proven by Auslander and
Goldman~\cite[Thm.\ 7.2]{AuGo60}. Since the class of
$A_{\xi}^{K}(x,y)^{\otimes\, n}=K\otimes_{R}A_{\xi}^{R}(x,y)^{\otimes\, n}$
is trivial in~$\Br (K)$, see~\cite[Expl.\ 15.3]{IalgKT} we
see that for regular integral domains~$R$ the classes of
the algebras~$A_{\xi}^{R}(x,y)$ lie in~$\Brtor{n}(R)$.

\medbreak

We will also use the notation~$A_{\xi}^{X}(x,y)$ for~$A_{\xi}^{R}(x,y)$,
where $X=\Spec R$. 
\end{emptythm}

\begin{emptythm}
\label{NormResSubSect}
{\it The norm residue homomorphism.}
Let~$n\geq 2$ be an integer and~$F$ a field
whose characteristic does not divide~$n$. We
denote by~$\MK{i}(F)$ the $i$-th Milnor $K$-group
of~$F$.

\smallbreak

As usual we denote ``pure'' symbols
in~$\MK{i}(F)$ by $\{ c_{1},\ldots ,c_{i}\}$,
where $c_{i}\in F^\times$.

\medbreak

By Hilbert~90 we get
from the long exact cohomology sequence associated
with the short exact sequence~(\ref{cdot-nEq}) an isomorphism
$$
R_{1,n}^{F}\, :\; F^{\times}/F^{\times\, n}\,
  \xrightarrow{\;\simeq\:}\,\HM^{1}(F,\roun_{n})\; ,\;\;
      \{ x\}\,\longmapsto\, (x)\, .
$$
Taking the cup product this map induces the so
called norm residue homomorphism
$$
R_{2,n}^{F}\, :\;\MK{2}(F)/n\cdot\MK{2}(F)\,
\too\,\HM^{2}(F,\roun_{n}^{\otimes\, 2})\; ,\;\;
\{ x,y\}\, +n\cdot\MK{2}(F)\,\longmapsto\, (x)\cup (y)\, ,
$$
which is by the Merkurjev-Suslin Theorem~\cite{MeSu82}
an isomorphism.
\end{emptythm}

\begin{emptythm}
\label{CyclicAlgExpl}
{\bfseries Example.}
Assume that~$F$ contains a primitive $n$-th
root of unity~$\xi$. Then we have isomorphisms
$$
\Brtor{n}(F)\,\xrightarrow{\;\simeq\:}\,
\HM^{2}(F,\roun_{n})\,\xrightarrow{\tau_{\xi}^{F}}\,
\HM^{2}(F,\roun_{n}^{\otimes\, 2})\, .
$$
The isomorphism~$\tau_{\xi}^{F}$ is induced by the $\Gamma_{F}$-module isomorphism
$\roun_{n}\xrightarrow{\simeq}\roun_{n}\otimes_{\Z}\roun_{n}$,
$\eta\,\longmapsto\,\eta\otimes\xi$, and so depends on the
choice of the primitive $n$-th root of unity~$\xi$. One can
show, see \eg~\cite[Prop.\ 4.7.1]{CSAaGC}, that the composition
of these isomorphisms maps the class of the algebra
$A_{\xi}^{F}(x,y)$ to~$(x)\cup (y)$.
\end{emptythm}

\begin{emptythm}
\label{mu-nNotation}
If the base field~$F$ contains a primitive $n$-th root of unity~$\xi$
we have also the isomorphism of \'etale sheaves
$\ru_{n}\xrightarrow{\simeq}\ru_{n}\otimes\ru_{n}$. It induces an
isomorphisms
$$
\tau^{X}_{\xi}\, :\;\HM^{2}_{\et}(X,\ru_{n})\,\xrightarrow{\;\simeq\;}\,
\HM^{2}_{\et}(X,\ru_{n}^{\otimes\, 2})
$$
for all $F$-schemes~$X$, which depends on the choice of the $n$-th root
of unity~$\xi$.
\end{emptythm}

\goodbreak
\section{On the Brauer group of the product of a torus
with certain regular and integral schemes}
\label{BrTorusSect}\bigbreak

\begin{emptythm}
\label{CycleComplSubSect}
{\it Cycle complexes.}
Let~$F$ be a field, and
$$
E\,\longmapsto\,\CM{\ast}(E)\,=\;\bigoplus\limits_{i\in\Z}\CM{i}(E)\, ,
$$
where $E\supseteq F$ is a field extension, a cycle module over
the field~$F$ in the sense of Rost~\cite{Ro96}. For instance,
Milnor $K$-theory $E\mapsto\MK{\ast}(E)$ or Milnor $K$-theory
modulo~$n$
$$
E\,\longmapsto\,\MK{\ast}(E)/n\cdot\MK{\ast}(E)\,
:=\;\bigoplus\limits_{i\geq 0}\MK{i}(E)/n\cdot\MK{i}(E)
$$
for some integer~$n\geq 2$ are examples of cycle modules.
For brevity we denote the later cycle module by~$\MK{\ast}/n$.

\smallbreak

If~$X$ is an integral scheme of finite type over~$F$ (or more general
an excellent scheme) Rost~\cite{Ro96} has constructed following Kato's
construction~\cite{Ka86} for Milnor $K$-theory a complex
$\CK^{\smb}(X,\CM{l})$:
$$
\xymatrix{
\qquad\CM{l}(F(X)) \ar[r]^-{d_{X}^{0}} &
{\bigoplus\limits_{x\in X^{(1)}}}\CM{l-1}(F(x)) \ar[r]^-{d_{X}^{1}} &
{\bigoplus\limits_{x\in X^{(2)}}}\CM{l-2}(F(x)) \ar[r] & \,\ldots\, ,
}
$$
for all integers~$l$, where~$X^{(i)}$ denotes the set of points of codimension~$i$
in~$X$, $F(x)$~the residue field of~$x\in X$, and  the
differentials~$d^{i}_{X}$ are induced by the second residue maps.
We consider~$\CK^{\smb}(X,\CM{l})$ as a cohomological complex
with $\bigoplus\limits_{x\in X^{(i)}}\CM{l-i}(F(x))$ in degree~$i$, and
denote the cohomology group in degree~$i$ by~$\HM^{i}(X,\CM{l})$.

\smallbreak

We will make use of the fact proven by Rost~\cite[Prop.\ 8.6]{Ro96}
that these groups are homotopy invariant, \ie the pull-back
$p^{\ast}:\HM^{i}(X,\CM{l})\too\HM^{i}(X\times_{F}\A^{1}_{F},\CM{l})$
along the projection $p:X\times_{F}\A^{1}_{F}\too X$ is an isomorphism
for all~$i\geq 0$.
\end{emptythm}

\begin{emptythm}
\label{cycleComplSubSect}
{\it The $\MK{2}/n$-unramified cohomology of a split torus.}
Let~$n\geq 2$ be an integer and $T=\Spec F[t_{1}^{\pm 1},\ldots ,t_{d}^{\pm 1}]$
a split torus of rank~$d\geq 1$. It follows from homotopy invariance
and the localization sequence that
\begin{multline}
\label{torusEq}
\HM^{0}(T,\MK{2}/n)\,\simeq\,
q_{T}^{\ast}\big(\MK{2}(F)/n\cdot\MK{2}(F)\,\big)\;\bigoplus\;\\
\left(\bigoplus\limits_{i=1}^{d}F^{\times}/F^{\times\, n}\,\cdot\{ t_{i}\}\right)
\;\bigoplus\;\left(\bigoplus\limits_{1\leq i<j\leq d}
\Z/\Z n\cdot\{ t_{i},t_{j}\}\right)\, ,
\end{multline}
where $q_{T}\colon T\too\Spec F$ is the structure morphism, and
$F^{\times}/F^{\times\, n}\cdot\{ t_{i}\}$ denotes subgroup
of $\MK{2}(F(T))/n\cdot \MK{2}(F(T))$ generated by
all $\{ a,t_{i}\}$ with $a\in F^\times$.
\end{emptythm}

\begin{emptythm}
\label{mainDiagSubSect}
{\it A commutative diagram.}
Let~$X$ be a regular and integral $F$-scheme of finite type and $n\geq 2$
be an integer which is not divisible by~$\khar F$. We set for brevity
$\mk{i}(E):=\MK{i}(E)/n\cdot\MK{i}(E)$ for all field extensions~$E\supseteq F$.

\smallbreak

We have then a commutative diagram, see Bloch~\cite[Thm.\ 2.3]{Bl80},
\begin{equation}
\label{mainDiag}
\xymatrix{
{\mk{2}(F(X))} \ar[r]^-{d_{X}^{0}} \ar[d]_-{R_{2,n}^{F(X)}}^-{\simeq} &
   {\bigoplus\limits_{x\in X^{(1)}}}\mk{1}(F(x)) \ar[r]^-{d_{X}^{1}}
       \ar[d]_-{(R_{1,n}^{F(x)})_{x\in X^{(1)}}}^-{\simeq} &
          {\bigoplus\limits_{x\in X^{(2)}}}\Z/\Z n \ar[d]^-{=}
\\
\HM^{2}(F(X),\roun_{n}^{\otimes\, 2}) \ar[r]_-{\delta_{X}^{0}} &
   {\bigoplus\limits_{x\in X^{(1)}}}\HM^{1}(F(x),\roun_{n}) \ar[r]_-{\delta_{X}^{1}} &
       {\bigoplus\limits_{x\in X^{(2)}}}\Z/\Z n \rlap{\, ,}
}
\end{equation}
where the lower row comes from the Bloch-Ogus resolution~\cite{BlOg74}.
By the Merkurjev-Suslin Theorem, see~\ref{NormResSubSect},
we know that~$R_{2,n}^{F(X)}$, which maps a symbol~$\{ x,y\}$
to the product~$(x)\cup (y)$, is an isomorphism. Hence
all vertical maps in this diagram are isomorphisms.

\smallbreak

Consequently, the norm residue homomorphism induces
an isomorphism between the kernel of $d_{X}^{0}:\mk{2}(F(X))
\too\bigoplus\limits_{x\in X^{(1)}}\mk{1}(F(x))$ and the
kernel of $\delta_{X}^{0}:\HM^{2}(F(X),\roun_{n}^{\otimes\, 2})
\too\bigoplus\limits_{x\in X^{(1)}}\HM^{1}(F(x),\roun_{n})$.

\medbreak

By the same diagram and~(\ref{torusEq}) we know that the kernel
of~$\delta^{0}_{T}$ is equal to the direct sum
$$
q_{T}^{\ast}\big(\HM^{2}(F,\roun_{n}^{\otimes\, 2})\big)\;\bigoplus\;
\left(\bigoplus\limits_{i=1}^{d}(F^{\times}/F^{\times\, n})\cup (t_{i})\right)
\;\bigoplus\;\left(\bigoplus\limits_{1\leq i<j\leq d}
\Z/\Z n\cdot (t_{i})\cup (t_{j})\right)\, ,
$$
where $(F^{\times}/F^{\times\, n})\cup (t_{i})$ means the subgroup
of~$\HM^{2}(F(T),\roun^{\otimes\, 2})$ generated by all cup
products~$(a)\cup (t_{i})$,
$a\in F^{\times}$, and $q_{T}:T\too\Spec F$ is the structure morphism.

\end{emptythm}

\begin{emptythm}
\label{SplitTorusSubSect}
{\it The unramified $\MK{2}/n$-cohomology of the product of a split
torus with certain regular schemes.}
We fix a field~$F$ and an integer~$n\geq 2$.
We assume first that the field~$F$ contains a primitive
$n$-th root of unity~$\xi$, and that the torus~$T$ is split, say
$T=\Spec F[t_{1}^{\pm 1},\ldots ,t_{d}^{\pm 1}]$, where $d=\rank T\geq 1$,
or that~$T$ is the ``trivial'' torus~$\Spec F$.

\smallbreak

Let further~$X$ be a regular and integral $F$-scheme of finite type
with structure morphism $q_{X}\colon X\too\Spec F$.
We make the following assumptions on~$X$:

\smallbreak

\begin{itemize}
\item[{\bf (H)}]
The natural homomorphism
$q_{T}^{\ast}:\MK{i}(F)/n\cdot\MK{i}(F)\too\HM^{0}(X,\MK{i}/n)$
is an isomorphism for~$i=1,2$.
\end{itemize}

\medbreak

\begin{emptythm}\label{SplitTorusSubSect2}
{\bfseries Lemma.}
{\it
Under these assumptions the pull-back homomorphism
$$
p_{T}^{\ast}\, :\;\HM^{0}(T,\MK{i}/n)\,\too\,\HM^{0}(X\times_{F}T,\MK{i}/n)
$$
along the projection $p_{T}\colon X\times_{F}T\too T$ is an isomorphism for~$i=1,2$.
Hence by the commutative diagram~(\ref{mainDiag})
also the pull-back homomorphism
$$
p_{T}^{\ast}\, :\;\Ker\delta^{0}_{T}\,\too\,\Ker\delta^{0}_{X\times_{F}T}\, .
$$
is an isomorphism.
}
\end{emptythm}
\begin{proof}
We prove this by induction on the rank of~$T$. The rank zero case,
\ie $T=\Spec F$, is assumption~{\bf (H)}. So let~$d=\rk T\geq 1$ and
$S:=\Spec F[t_{1}^{\pm 1},\ldots ,t_{d-1}^{\pm 1}]$. Then we get
from the localization sequences for the open embeddings
$\iota:T\hookrightarrow S[t_{d}]$ and $\id_{X}\times\iota :X\times_{F}T
\hookrightarrow X\times_{F}S[t_{d}]$ a commutative diagram whose
rows are exact
$$
\xymatrix{
\HM^{0}(S[t_{d}],\MK{i}/n)\;\; \ar@{>->}[r]^-{\iota^{\ast}} \ar[d]_-{(p_{S}[t_{d}])^{\ast}} &
    \HM^{0}(T,\MK{i}/n) \ar@{->>}[r]^-{\partial_{T}} \ar[d]^-{p_{T}^{\ast}} &
        \HM^{0}(S,\MK{i-1}/n) \ar[d]^-{p_{S}^{\ast}}
\\
\HM^{0}(X\times_{F}S[t_{d}],\MK{i}/n)\;\; \ar@{>->}[r]^-{(\id_{X}\times\iota)^{\ast}} &
   \;\HM^{0}(X\times_{F}T,\MK{i}/n) \ar@{->>}[r]^-{\partial_{X\times T}} &
      \HM^{0}(X\times_{F}S,\MK{i-1}/n)\, ,
}
$$
where $p_{S}:X\times_{F}S\too S$ is the projection. As indicated the left
horizontal arrows~$\iota^{\ast}$ and $(\id_{X}\times\iota)^{\ast}$ are injective
and both connecting homomorphisms~$\partial_{T}$ and~$\partial_{X\times T}$ are
surjective. In fact these are split surjections. The splitting is induced
by multiplication with the symbol of length one~$\{ t_{d}\}$.

\smallbreak

The right hand vertical
arrow $p_{S}^{\ast}$ is an isomorphism for~$i=1$ since~$X$ and~$X\times_{F}S$
are irreducible and so $\HM^{0}(S,\MK{0}/n)=\HM^{0}(X\times_{F}S,\MK{0}/n)=\Z/\Z n$,
and for~$i=2$ by induction. The left hand vertical arrow~$(p_{S}[t_{d}])^{\ast}$ is an
isomorphism by induction and homotopy invariance, and so it follows
by the snake lemma that also the vertical arrow~$p_{T}^{\ast}$ in the middle
of the diagram is an isomorphism. We are done.
\end{proof}

\medbreak

We use this to describe~$\Brtor{n}(T)$ and $\HMtor{n}^{2}_{\et}(X\times_{F}T,\Gm)$.
\end{emptythm}

\begin{emptythm}
\label{BrSplitTorusSubSect}
{\it The Brauer group of a split torus~$T$.}
Let~$T$ and~$F$ be as above. Since~$\Pic T=0$ as~$T$ is a split torus
we know (\cf~\ref{BrauerGrSubSect}) that $\Brtor{n}(T)=\HM^{2}_{\et}(T,\ru_{n})$,
and so by Auslander and Goldman~\cite[Thm.\ 7.2]{AuGo60} the pull-back homomorphism
$$
\iota^{\ast}_{T}\, :\;\HM^{2}_{\et}(T,\ru_{n})\,\too\,
\HM^{2}_{\et}(F(T),\ru_{n})=\HM^{2}(F(T),\roun_{n})
$$
is injective, where $\iota_{T}:\Spec F(T)\too T$ denotes the generic point. By the very definition
of the Bloch-Ogus complex the image of $\iota^{\ast}_{T}:\HM^{2}_{\et}(T,\ru_{n}^{\otimes\, 2})
\too\HM^{2}(F(T),\ru_{n}^{\otimes\, 2})$ lies in the kernel of~$\delta^{0}_{T}$.
Hence we have a commutative diagram where all vertical arrows are monomorphisms
\begin{equation}
\label{T-Diag}
\xymatrix{
\HM^{2}_{\et}(F,\roun_{n}) \ar[d]^-{q_{T}^{\ast}} \ar[rr]^-{\tau^{F}_{\xi}} & & \HM^{2}(F,\roun_{n}^{\otimes\, 2})
    \ar[d]^-{q_{T}^{\ast}} &
\\
\HM^{2}_{\et}(T,\ru_{n}) \ar[rr]^-{\tau^{T}_{\xi}}_-{\simeq} \ar[d]^-{\iota^{\ast}_{T}} & &
   \HM^{2}_{\et}(T,\ru_{n}^{\otimes\, 2}) \ar[r]^-{\iota^{\ast}_{T}} \ar[d]^-{\iota^{\ast}_{T}} &
             \Ker\delta^{0}_{T} \ar[ld]^-{\subseteq}
\\
\HM^{2}(F(T),\roun_{n}) \ar[rr]^{\tau^{F(T)}_{\xi}}_-{\simeq} & & \HM^{2}(F(T),\roun_{n}^{\otimes\, 2})\, . &
}
\end{equation}
Therefore $\iota_{T}^{\ast}\circ\tau^{T}_{\xi}:\HM^{2}_{\et}(T,\ru_{n})\too\Ker\delta^{0}_{T}$
is injective. We claim that it is also surjective.

\smallbreak

In fact, let $x,y$ be units in~$F[T]^{\times}$ then the class of the Azumaya algebra
$A^{T}_{\xi}(x,y)$ is in~$\Brtor{n}(X)=\HM^{2}_{\et}(X,\ru_{n})$, see~\ref{AzExplSubSect}.
By Example~\ref{CyclicAlgExpl} we have
$$
\tau^{F(T)}_{\xi}\big(\,\iota_{T}^{\ast}([A^{T}_{\xi}(x,y)]\,\big)\, =\, (x)\cup (y)\, .
$$
Since as seen above~$\Ker\delta^{0}_{T}$ is the direct sum of
$q_{T}^{\ast}\big(\HM^{2}(F,\roun_{2}^{\otimes\, 2})\big)$ and the subgroup generated
by all $(a)\cup (t_{i})$ and $(t_{i})\cup (t_{j})$ with $a\in F^{\times}$ and $1\leq i,j\leq d$
we get therefore from diagram~(\ref{T-Diag}) our claim.

\medbreak

We have proven the following theorem which is due to
Magid~\cite{Ma78} if~$F$ is an algebraically closed
field of characteristic~$0$, see also P.~Gille and
A.~Pianzola~\cite[Sect.\ 4]{GiPa08} for another proof
in this case.
\end{emptythm}

\begin{emptythm}
\label{torusThm}
{\bfseries Theorem.}
{\it
Let~$n\geq 2$ be an integer, $T=\Spec F[t_{1}^{\pm 1},\ldots ,t_{d}^{\pm 1}]$
be a split torus of rank~$d\geq 1$ over a field~$F$ which contains a
primitive $n$-th root of unity~$\xi$. Then $\Brtor{n}(T)$ is isomorphic to
$$
\Brtor{n}(F)\oplus\;\bigoplus\limits_{i=1}^{d}\sum\limits_{a\in F^{\times}}\Z/\Z n\cdot [A_{\xi}^{T}(a,t_{i})]
$$
if~$d=1$ and to
$$
\Brtor{n}(F)\oplus B(T)\,\oplus\;
      \bigoplus\limits_{i=1}^{d}\sum\limits_{a\in F^{\times}}\Z/\Z n\cdot [A_{\xi}^{T}(a,t_{i})]
           \;\oplus\;\bigoplus\limits_{1\leq i<j\leq d}\Z/\Z n\cdot
        [A^{T}_{\xi}(t_{i},t_{j})]
$$
if~$d=2$. In particular, if $F^{\times}=F^{\times\, n}$ we have a natural isomorphism
$\Brtor{n}(F)\xrightarrow{\simeq}\Brtor{n}(T)$ if the rank of~$T$ is~$1$,
and an isomorphism
$$
\Brtor{n}(T)\,\simeq\,\Brtor{n}(F)\oplus\,\bigoplus\limits_{1\leq i<j\leq d}
\Z/\Z n\cdot [A^{T}_{\xi}(t_{i},t_{j})]
$$
if the rank of~$T$ is~$d\geq 2$.
}
\end{emptythm}

\begin{emptythm}
\label{product-TorusSubSect}
Let (as above)~$n\geq 2$ be an integer and~$F$ a field which
contains a primitive $n$-th root of unity. Let further~$T$ be
a split torus or equal~$\Spec F$ and~$X$ a regular and integral
$F$-scheme of finite type which satisfies the condition~{\bf (H)}
in~\ref{SplitTorusSubSect} and also
\begin{itemize}
\item[{\bf (nd)}]
$\Pic X$ is $n$-divisible.
\end{itemize}

\medbreak

\noindent
We set $Y=X\times_{F}T$ and denote $p_{T}:Y\too T$ the projection
(which is the structure morphism if $T=\Spec F$). Then as~$\Pic X$
is $n$-divisible it follows by homotopy invariance and localization
sequence that also~$\Pic Y$ is $n$-divisible. Therefore
$\HM^{2}_{\et}(Y,\ru_{n})=\HMtor{n}^{2}_{\et}(Y,\Gm)$ and so
by Grothendieck~\cite[Chap.\ II, Cor.\ 1.8]{Gr68}, see
also~\cite[Chap.\ III, Expl.\ 2.22]{EC}, we know that the pull-back
$\iota^{\ast}_{Y}:\HM^{2}_{\et}(Y,\ru_{n})\too\HM^{2}(F(Y),\roun_{n})$
along the generic point $\iota_{Y}:\Spec F(Y)\too Y$ is a monomorphism.

\smallbreak

We have now the following commutative diagram
$$
\xymatrix{
\HM^{2}_{\et}(T,\ru_{n}) \ar[rr]^-{\tau^{T}_{\xi}}_-{\simeq} \ar[d]_-{p_{T}^{\ast}} & &
   \HM^{2}(T,\ru_{n}^{\otimes\, 2}) \ar[r]^-{\iota_{T}^{\ast}}_-{\simeq} \ar[d]_-{p_{T}^{\ast}} &
       \Ker\delta^{0}_{T} \ar[d]^-{p_{T}^{\ast}}
\\
\HM^{2}_{\et}(Y,\ru_{n}) \ar[rr]^-{\tau^{Y}_{\xi}}_-{\simeq} \ar[d]_-{\iota_{Y}^{\ast}} & &
   \HM^{2}(Y,\ru_{n}^{\otimes\, 2}) \ar[r]^-{\iota_{Y}^{\ast}} \ar[d]_-{\iota_{Y}^{\ast}} &
       \Ker\delta^{0}_{Y} \ar[ld]^-{\subseteq}
\\
\HM^{2}(F(Y),\roun_{n}) \ar[rr]^-{\tau_{\xi}^{F(Y)}}_-{\simeq} & &
    \HM^{2}(F(Y),\roun_{n}^{\otimes\, 2})\, . &
}
$$
The vertical arrow on the right hand side
$p_{T}^{\ast}:\Ker\delta^{0}_{T}\too\Ker\delta^{0}_{Y}$ is an
isomorphism by the lemma~\ref{SplitTorusSubSect2} above,
and we have shown in~\ref{BrSplitTorusSubSect} that also~$\iota_{T}^{\ast}$
is an isomorphism. Hence since $\iota_{Y}^{\ast}:\HM^{2}_{\et}(Y,\ru_{n})
\too\HM^{2}(F(Y),\roun_{n})$ is injective we get from this diagram
that also $p_{T}^{\ast}:\HM^{2}_{\et}(T,\ru_{n})\too\HM^{2}_{\et}(Y,\ru_{n})$
is an isomorphism.

\medbreak

We have proven:

\medbreak

\begin{emptythm}
{\bfseries Theorem.}
{\it
Let~$n$ be an integer~$\geq 2$, $F$~a field which contains a primitive
$n$-th root of unity and~$T$ either a split torus or equal $\Spec F$. Let further~$X$
be a regular and integral $F$-scheme, such that~(i) $\Pic X$ is $n$-divisible,
and~(ii) the natural homomorphism $\MK{i}(F)/n\cdot\MK{i}(F)\too\HM^{0}(X,\MK{i}/n)$
is an isomorphism for $i=1,2$. Then the pull-back homomorphism
$$
p_{T}^{\ast}\, :\;\Brtor{n}(T)\,\too\,\HMtor{n}^{2}_{\et}(X\times_{F}T,\Gm)
$$
where $p_{T}:X\times_{F}T\too T$ is the projection, is an isomorphisms.
}
\end{emptythm}

\smallbreak

Note that if $\Br (F)\too\HM^{2}_{\et}(X,\Gm)$ is injective (\eg if
$X(F)\not=\emptyset$) and the natural homomorphism
$\Brtor{n}(F)=\HM^{2}(F,\roun_{n})\too\HM^{2}_{\et}(X,\ru_{n})$ is
an isomorphism then it follows from the long exact cohomology
sequence for the exact sequence~(\ref{cdot-nEq}) that~$\Pic X$
is $n$-divisible. Hence in this case condition~{\bf (nd)} 
is a necessary condition.

\medbreak

We consider now the case that the torus~$T$ is not split. Here we show
the following generalization.
\end{emptythm}

\begin{emptythm}
\label{mainThm}
{\bfseries Theorem.}
{\it
Let~$F$ be a field and~$n\geq 2$ an integer which is not divisible
by~$\khar F$. Let further~$X$ be a regular and geometrically integral
$F$-scheme of finite type, and~$T$ a $F$-scheme, such that either~$T_{s}$
is isomorphic as $F_{s}$-scheme to a torus, or~$T$ is the ``trivial''
torus~$\Spec F$. We assume that
\begin{itemize}
\item[(i)]
$\Pic X_{s}=0$,

\smallbreak

\item[(ii)]
$\HM^{0}(X_{s},\MK{i}/n)=0$ for~$i=1,2$, and

\smallbreak

\item[(iii)]
$F_{s}[X_{s}]^{\times}=F_{s}^{\times}$.
\end{itemize}
Then the natural homomorphism
$$
p_{T}\, :\;\Brtor{n}(T)\,\too\,\HMtor{n}^{2}_{\et}(X\times_{F}T,\Gm)
$$
along the projection $p_{T}:X\times_{F}T\too T$ is an isomorphism.
}

\begin{proof}
We show first that $\Brtor{n}(F)\too\Brtor{n}(X)$ is an
isomorphism. If~$F$ is separably closed this is proven
in~\ref{product-TorusSubSect} above. In the general
case we get since~$\Pic X_{s}=0$ from the Hochschild-Serre
spectral sequence, see~\ref{HSSpSeqSubSect}, an exact sequence
$$
\xymatrix{
0 \ar[r] & \HM^{2}(F,F_{s}[X_{s}]^{\times}) \ar[r] &
   \HM^{2}_{\et}(X,\Gm) \ar[r] & \HM^{2}_{\et}(X_{s},\Gm)^{\Gamma_{F}}}\, ,
$$
where~$\Gamma_{F}$ is the absolute Galois group of~$F$.
By case of a separably closed base field we know
that $\HMtor{n}^{2}_{\et}(X_{s},\Gm)=0$, and by our assumption~(iii)
we have $F_{s}[X_{s}]^{\times}=F_{s}^{\times}$, and so
$\HM^{2}(F,F_{s}[X_{s}]^{\times})=\Br (F)$. Hence we get
from the exact sequence above the claimed isomorphism
$\Brtor{n}(F)\xrightarrow{\simeq}\HMtor{n}^{2}_{\et}(X,\Gm)$.

\medbreak

Assume now that~$T$ is not the ``trivial'' torus. Then
if~$F=F_{s}$ is separably closed this is the assertion
of the theorem in~\ref{product-TorusSubSect}. In the general
case, we know by homotopy invariance and the localization sequence
that $\Pic T_{s}=0$ and $\Pic (X\times_{F}T)_{s}=0$. Hence we get from
the Hochschild-Serre spectral sequence, see~\ref{HSSpSeqSubSect},
a commutative diagram whose rows are exact
$$
\xymatrix{
\HM^{2}(F,F_{s}[T_{s}]^{\times})\;\; \ar@{>->}[r]^-{\alpha} \ar[d]^-{p^{\ast}_{T_{s}}}_{\simeq} &
 \Br (T) \ar[r]^-{\beta} \ar[d]^-{p^{\ast}_{T}} &  \Br (T_{s})^{\Gamma_{F}} \ar[r]^-{\gamma}
     \ar[d]^-{p^{\ast}_{T_{s}}} & \HM^{3}(F,F_{s}[T_{s}]^{\times}) \ar[d]^-{p^{\ast}_{T_{s}}}_{\simeq}
\\
\HM^{2}(F,F_{s}[Y_{s}]^{\times})\;\; \ar@{>->}[r]^-{\alpha'} & \HM^{2}_{\et}(Y,\Gm) \ar[r]^-{\beta'} &
  \HM^{2}_{\et}(Y_{s},\Gm)^{\Gamma_{F}} \ar[r]^-{\gamma'} & \HM^{3}(F,F_{s}[Y_{s}]^{\times})\, ,
}
$$
where we have set~$Y=X\times_{F}T$ and~$p_{T}:X\times_{F}T\too T$
and~$p_{T_{s}}:X_{s}\times_{F_{s}}T_{s}\too T_{s}$
are the respective projections. In particular, the arrows~$\alpha$
and~$\alpha'$ on the left hand side are injective.
As indicated the utmost left and right vertical arrows 
are isomorphisms, since by Rosenlicht~\cite[Thms.\  1 and 2]{Ros61},
see also~\cite[Lem.\ 10]{CoThSa77}, we know that the morphism~$p_{T_{s}}$
induces an isomorphism $F_{s}[T_{s}]^{\times}\xrightarrow{\simeq}
F_{s}[Y_{s}]^{\times}$. By the split case $p_{T}^{\ast}:\Br (T_{s})\too\HM^{2}_{\et}(Y_{s},\Gm)$
induces an isomorphism between the subgroups of elements of order dividing~$n$.
This implies by a straightforward diagram chase that
$p_{T}^{\ast}:\Brtor{n}(T)\too\HMtor{n}^{2}_{\et}(Y,\Gm)$
is injective.

\smallbreak

To show that this map is also surjective, let~$x'\in\HMtor{n}^{2}_{\et}(Y,\Gm)$.
Then by the split case there is $y\in\Brtor{n}(T_{s})^{\Gamma_{F}}$,
such that $p_{T}^{\ast}(y)=\beta' (x')$. We have $p_{T}^{\ast}(\gamma (y))=
\gamma'(p_{T}^{\ast}(y))=\gamma'(\beta'(x'))=0$ and so there exists a~$z\in\Br(T)$,
such that $\beta (z)=y$.

\smallbreak

It follows then that
$\beta'(x'-p_{T}^{\ast}(z))=0$ and so we get
$x'=p_{T}^{\ast}(z+\alpha (w))$ for some
$w\in\HM^{2}(F,F_{s}[T_{s}]^{\times})$. Since as
observed above $p_{T}^{\ast}:\Brtor{n}(T)\too\HMtor{n}^{2}_{\et}(Y,\Gm)$
is a monomorphism the assumption $n\cdot x'=0$
implies also $n\cdot (z+\alpha (w))=0$ and so the
homomorphism $p_{T}^{\ast}:\Brtor{n}(T)\too\HMtor{n}^{2}_{\et}(Y,\Gm)$
is also surjective, hence an isomorphism.
\end {proof}

\medbreak

\begin{emptythm}
{\bfseries Corollary.}
{\it
Let~$n$ and~$F$ be as above, and~$Z$ a regular
and integral $F$-scheme which contains
an open $F$-subscheme~$X$ of finite type which satisfies conditions
(i)--(iii) in the theorem above, and let~$T$ be a $F$-scheme,
such that either~$T_{s}$ is isomorphic to a torus as $F_{s}$-scheme,
or $T=\Spec F$. Then $\Brtor{n}(T)\too\HMtor{n}^{2}_{\et}(Z\times_{F}T,\Gm)$
is an isomorphism.
}
\end{emptythm}
\begin{proof}
Since~$Z$ and so also $Z\times_{F}T$ are regular and integral
schemes the natural homomorphisms
$\HMtor{n}^{2}_{\et}(Z,\Gm)\too\HMtor{n}^{2}_{\et}(X,\Gm)$
and $\HMtor{n}^{2}_{\et}(Z\times_{F}T,\Gm)\too
\HMtor{n}^{2}_{\et}(X\times_{F}T,\Gm)$ are injective
by Grothendieck~\cite[Chap.\ II, Cor.\ 1.8]{Gr68}.
\end{proof}
\end{emptythm}

\begin{emptythm}
\label{p-Rem}
{\bfseries Remark.}
If $\khar F=p>0$ and~$X$ is a regular affine and geometrically
integral $F$-scheme of dimension~$\geq 2$ then the cokernel
of the homomorphism $\Br (F)\too\Br (X)$ contains
a non empty $p$-torsion and $p$-divisible subgroup, and so
is in particular an infinite abelian group, see our
Theorem~\ref{mainp-Thm}.
\end{emptythm}

\begin{emptythm}
\label{GeeometricCriterionSubSect}
{\it A geometric criterion.}
We give here a geometric condition which implies~{\bf (H)}
in~\ref{SplitTorusSubSect}. To formulate it we introduce
the following notation. Let~$E$ be a field and~$\A^{1}_{E}$
the affine line over~$E$. We denote for~$i=0,1$ by
$s_{i}:\Spec E\too\A^{1}_{E}$ the points~$0$ and~$1$, respectively.

\medbreak

Let ~$F$ be a field and~$X$ a regular and integral
$F$-scheme of finite type. Let~$\iota_{X}:\Spec F(X)\too X$ be the
generic point and $q_{X}:X\too\Spec F$ the structure morphism.
 Assume that there is an $F$-rational point
$x:\Spec F\too X$, and there are $F$-morphisms
$h_{j}:\A^{1}_{F(X)}\too X$, $j=1,\ldots ,l$, such that~(a)
$h_{1}\circ s_{0}$ factors through~$x$, (b)
$h_{j}\circ s_{1}=h_{j+1}\circ s_{0}$ for
$j=1,\ldots ,l-1$, and~(c) $h_{l}\circ s_{1}=\iota_{X}$.
Then
$$
q_{X}^{\ast}\, :\;\CM{k}(F)\,\too\,\HM^{0}(X,\CM{k})
$$
is an isomorphism for all~$k\in\Z$ and all cycle
modules~$\CM{\ast}$ over~$F$.

\smallbreak

In fact, let $q_{\A^{1}_{F(X)}}:\A^{1}_{F(X)}\too\Spec F(X)$
be the structure morphism. Then by homotopy invariance
the pull-back $q_{\A^{1}_{F(X)}}^{\ast}:\CM{k}(F(X))\too\HM^{0}(\A^{1}_{F(X)},\CM{k})$
is an isomorphism. Since $q_{\A^{1}_{F(X)}}\circ s_{i}=\id_{\Spec F(X)}$
for $i=0,1$ we have $s_{0}^{\ast}=q_{\A^{1}_{F(X)}}^{\ast}=s_{1}^{\ast}$.
Therefore $(h_{j}\circ s_{0})^{\ast}=(h_{j}\circ s_{1})^{\ast}$ for
all $j=1,\ldots ,l$, and so it follows from~(b) that
$$
(h_{1}\circ s_{0})^{\ast}\, =\, (h_{l}\circ s_{1})^{\ast}\, :\;
   \CM{k}(F(X))\,\too\,\HM^{0}(X,\CM{k})\, .
$$
But by~(c) we have $h_{l}\circ s_{1}=\iota_{X}$ and so
by the very definition of the pull-back homomorphism on
cycle complexes $(h_{l}\circ s_{1})^{\ast}$ and hence
also $(h_{1}\circ s_{0})^{\ast}$ is injective. This implies
by~(a) that $x^{\ast}:\HM^{0}(X,\CM{k})\too\CM{k}(F)$ is
injective and so $q_{X}^{\ast}:\CM{k}(F)\too\HM^{0}(X,\CM{k})$
an isomorphism since $q_{X}\circ x=\id_{\Spec F}$.
\end{emptythm}

\goodbreak
\section{Applications: Brauer groups of affine quadrics and
some reductive algebraic groups}
\label{ExplSect}\bigbreak

\begin{emptythm}
\label{affQuadSubSect}
{\it The Brauer group of an affine quadric.}
We assume here that~$F$ is a field of characteristic~$\not= 2$.
As above~$n$ is an integer which is not divisible by~$\khar F$.

\smallbreak

Let
$$
q\, =\,\sum\limits_{i=0}^{m}a_{i}x_{i}^{2}\, ,\; a_i\in F^\times\, ,\; m\geq 4\, ,
$$
be a regular quadratic form over the field~$F$. Denote
by~$X_{q}\hookrightarrow\P^{m}_{F}$
the corresponding projective quadric, and by~$X_{q,\aff}\subset X_{q}$
the open affine subscheme defined by $x_{0}\not=0$. This
is (isomorphic to) a non singular affine quadric
in~$\A^{m}_{F}$ given by the equation $a_{0}+\sum\limits_{j=1}^{m}a_{j}x_{j}^{2}=0$,
and so in particular a smooth and geometrically integral affine scheme.
We want to apply Theorem~\ref{mainThm} to compute the Brauer group of
this $F$-scheme, \ie we have to verify conditions~(i)--(iii)
there. Since this are assumptions on $F_{s}\times_{F}X_{q,\aff}$
we can assume for ease of notation that~$F=F_{s}$
is separably closed. Then~$q$ is a split quadratic form and
so we can assume that
$q=a_{0}x_{0}+x_{1}\cdot x_{2}+\sum\limits_{j=3}^{m}a_{j}x_{j}^{2}$.
The variety~$X_{q,\aff}$  is then isomorphic to the
affine quadric in~$\A^{m}_{F}$ given by the equation
$a_{0}+x_{1}\cdot x_{2}+\sum\limits_{j=3}^{m}a_{i}x_{i}^{2}=0$.

\smallbreak

Consider the open subscheme $x_{1}\not= 0$ of~$X_{q,\aff}$. It is
isomorphic to the spectrum of $F[x_{1}^{\pm 1},x_{3},\ldots ,x_{m}]$ and so every unit
of this open subscheme is equal $a\cdot x_{1}^{r}$ for some~$a\in F^{\times}$
and some~$r\in\Z$. Since~$x_{1}$ is not a unit in~$F[X_{q,\aff}]$
we see that $F[X_{q,\aff}]^{\times}=F^{\times}$, hence~(iii). To
show~(ii), \ie $\Pic X_{q,\aff}=0$, we recall that
$\CH_{m-1}(X_{q})\simeq\Pic X_{q}$ is generated by the class
of the hyperplane section $x_{0}=0$. Therefore by the localization
sequence we have $\Pic X_{q,\aff}\simeq\CH_{m-1}(X_{q,\aff})=0$.

\smallbreak

Finally we have to check~(ii), \ie $\HM^{2}(X_{q,\aff},\MK{i}/n)=0$
for $i=1,2$. For this we note first that the closed subscheme  in~$X_{q}$
defined by~$x_{0}=0$ is isomorphic to the projective quadric
$X_{q'}\hookrightarrow\P^{m-1}_{F}$, where
$q'=x_{1}\cdot x_{2}+\sum\limits_{j=3}^{m}a_{i}x_{i}^{2}$. Hence we have
an exact sequence of complexes
\begin{equation}
\label{CKEq}
\xymatrix{
\CK^{\smb}(X_{q'},\MK{i-1}/n)[-1]\;\; \ar@{>->}[r] &
     \CK^{\smb}(X_{q},\MK{i}/n) \ar@{->>}[r] & \CK^{\smb}(X_{q,\aff},\MK{i}/n)\, ,
}
\end{equation}
for all $i\geq 1$, where we denote for a complex~$K^{\smb}$ by~$K^{\smb}[-1]$
the shifted by~$-1$ complex: $(K^{\smb}[-1])^{j}=K^{j-1}$.

\smallbreak

To analyze the associated cohomology sequence for $i=1$ and~$2$
we recall that for a projective quadric~$Q$ of (Krull-)dimension~$\geq 1$
over the separably closed field~$F$ we have $\HM^{0}(Q,\MK{j}/l)=0$
for all $j\geq 1$ and~$l\geq 1$. In fact, the motive of~$Q$ is
isomorphic to a direct sum of Tate motives
$\bigoplus\limits_{k=0}^{\dim X_{q}}\ul{\Z}(k)^{\oplus\, r_{k}}$,
where the integers~$r_{k}$ are~$1$, except if~$\dim X_{q}$ is even
and~$k=\frac{\dim X_{q}}{2}$, in which case it is~$2$. Hence we have
$$
\HM^{0}(Q,\MK{j}/l)\, =\,\MK{j}(F)/l\cdot\MK{j}(F)
$$
for all~$j$, and the group $\MK{j}(F)/l\cdot\MK{j}(F)$
is zero if~$l\geq 1$ since~$F$ is separably closed.

\smallbreak

Hence we get from the long exact cohomology sequence associated
with~(\ref{CKEq}) a exact sequence
$$
\xymatrix{
0 \ar[r] & \HM^{0}(X_{q,\aff},\MK{i}/n) \ar[r]^-{\partial} &
    \HM^{0}(X_{q'},\MK{i-1}/n) \qquad\qquad\qquad\qquad\qquad
}
$$
$$
\xymatrix{
\qquad\qquad\qquad\qquad\qquad\qquad \ar[r] & \HM^{1}(X_{q},\MK{i}/n) \ar[r] &
           \HM^{1}(X_{q,\aff},\MK{i}/n)
}
$$
for all~$i\geq 1$. Since as observed above~$\HM^{0}(X_{q'},\MK{1}/n)=0$
we get form this exact sequence immediately $\HM^{0}(X_{q,\aff},\MK{2}/n)=0$.

\smallbreak

To show that also $\HM^{0}(X_{q,\aff},\MK{1}/n)=0$ we
note that $\HM^{i}(Y,\MK{i}/n)\simeq\Z/\Z n\otimes_{\Z}\CH^{i}(Y)$
for a smooth and integral $F$-scheme~$Y$ and so by what we have shown
above this exact sequence of cohomology groups gives a short
exact sequence
$$
\xymatrix{
0 \ar[r] & \HM^{0}(X_{q,\aff},\MK{1}/n) \ar[r] & \Z/\Z n \ar[r] &
     \Z/\Z n \ar[r] & 0\, ,
}
$$
which in turn implies that $\HM^{0}(X_{q,\aff},\MK{1}/n)$ is trivial.

\medbreak

This proves by Theorem~\ref{mainThm} and its corollary the following result.
\end{emptythm}

\begin{emptythm}
\label{affQuadThm}
{\bfseries Theorem.}
{\it
Let~$F$ be a field of characteristic~$\not= 2$, and $X_{q}\hookrightarrow\P^{m}_{F}$
be the projective quadric defined by the equation
$q=\sum\limits_{i=0}^{m}a_{i}x_{i}^{2}$
with $m\geq 4$ and $a_{i}\in F^{\times}$ for all $i=0,1,\ldots m$.
Let further~$X_{q,\aff}$ be the open affine quadric defined
by~$x_{0}\not= 0$ and $n\geq 2$ an integer which is not divisible by
the characteristic of~$F$. Then the natural homomorphisms
$$
\Brtor{n}(F)\,\too\,\Brtor{n}(X_{q,\aff})\quad\mbox{and}\quad
\Brtor{n}(F)\,\too\,\HMtor{n}^{2}_{\et}(X_{q},\Gm)
$$
are isomorphisms, and if~$X_{q,\aff}(F)\not=\emptyset$ also
the pull-back homomorphisms
$$
\Brtor{n}(T)\,\too\,\Brtor{n}(X_{q,\aff}\times_{F}T)\quad
\mbox{and}\quad\Brtor{n}(T)\,\too\,
\HMtor{n}^{2}_{\et}(X_{q},\Gm)
$$
along the respective projections to~$T$ are isomorphisms for
any $F$-scheme~$T$, such that~$T_{s}$ is isomorphic as $F_{s}$-scheme
to a torus.
}
\end{emptythm}
\medbreak

\begin{emptythm}
{\bfseries Remark.}
Merkurjev and Tignol~\cite[Thm.\ B]{MeTi95} have shown that for
any projective homogeneous variety~$X$ for a semisimple algebraic
group over a field~$F$ the natural homomorphism
$\Br (F)\too\HM^{2}_{\et}(X,\Gm)$ is surjective, and moreover they
have given a description of the kernel. In particular they proved,
see~\cite[Cor.\ 2]{MeTi95}, that for a projective quadric~$X_{q}$
over a field~$F$ as above~$\Br (F)\simeq\HM^{2}_{\et}(X_{q},\Gm)$.
\end{emptythm}

\begin{emptythm}
\label{scGrSubSect}
{\it Simply connected algebraic groups.}
Let~$G$ be a simply connected algebraic group
over the field~$F$. We assume that~$\khar F$ does not
divide an integer~$n\geq 2$. Then we have
$$
\HM^{0}(G_{s},\MK{i}/n)\, =\,\MK{i}(F_{s})/n\cdot\MK{i}(F_{s})\, =0
$$
for $i=1,2$, and also $\Pic G_{s}=0$, see \eg~\cite[Thm.\ 1.5]{Gi09}
for a proof. Moreover by Rosenlicht~\cite{Ros61} we know
that $F_{s}[G_{s}]^{\times}=F_{s}^{\times}$, and so
we conclude from Theorem~\ref{mainThm} and its corollary
the following.

\medbreak

\begin{emptythm}
\label{cor34}
{\bfseries Corollary.}
{\it
Let~$G$ be an algebraic variety (not necessarily an algebraic group)
over a field~$F$ such that the variety $G_{s}$ is isomorphic to
a simply connected algebraic group over~$F_{s}$,
and~$n\geq 2$ an integer which is not divisible by the characteristic
of~$F$. Then the pull-back homomorphism
$$
p_{T}^{\ast}\, :\;\Brtor{n}(T)\,\too\,\Brtor{n}(G\times_{F}T)
$$
along the projection $p_{T}:G\times_{F}T\too T$ is an isomorphism
in the following cases:
\begin{itemize}
\item[(i)]
$T=\Spec F$ and so~$p_{T}$ is the structure morphism $G\too\Spec F$; or

\smallbreak

\item[(ii)]
the $F$-scheme $T$~is a form of a torus, \ie $F_{s}\times_{F}T$
is isomorphic as $F_{s}$-scheme to a torus over~$F_{s}$, and $G(F)\ne\emptyset$.
\end{itemize}

\smallbreak

\noindent
In particular, if $\khar F=0$ we have an isomorphism
$$
p_{T}^{\ast}\, :\;\Br (T)\,\xrightarrow{\;\simeq\;}\,\Br (G\times_{F}T)
$$
for such~$G$ and~$T$.

Moreover, if $\overline G$ is a smooth
compactification of $G$, then the natural pull-back
homomorphism $\Brtor{n}(F)\,\too\,\HMtor{n}^{2}_{\et}(\ol{G},\Gm)$
is an isomorphism.
}
\end{emptythm}
\end{emptythm}

\begin{emptythm}
\label{semisimpleRem}
{\bfseries Remark.}
The first case $T=\Spec F$ was proven by Iversen~\cite{Iv76}
using topological methods if~$F$ is an algebraically closed field of
characteristic~$0$, and by the first author in~\cite{Gi09} for arbitrary fields.
\end{emptythm}

\begin{emptythm}
\label{semisimpleGrSubSect}
{\it Product of a torus with an arbitrary semisimple
algebraic group.}
Let (as above)~$F$ be a field and~$n\geq 2$ an integer
which is not divisible by~$\khar F$.

\smallbreak

Let~$G$ be an arbitrary connected semisimple algebraic group over~$F$
and~$T$ an $F$-scheme, such that~$T_{s}$ is isomorphic to a
torus over~$F_{s}$.
Let $\pi\colon G_{\sico}\too G$ be the
simply connected cover, see Tits~\cite[Prop.\ 2]{Ti66}.
We denote by $\Pi_{G}$ the fundamental group of $G$, \ie more
precisely the fundamental group of $F_{s}\times_{F}G$.
We assume that $\khar F$ does not divide the order of~$\Pi_{G}$.
Then the induced morphism $F_{s}\times_{F}G_{\sico}\too F_{s}\times_{F}G$ is
a Galois cover with group~$\Pi_{G}$. In particular, the
function field extension $F_{s}(G_{\sico})\supseteq F_{s}(G)$
is a Galois extension of degree~$|\Pi_{G}|$. Since~$F_{s}\times_{F}G$
is connected this implies that~$F(G_{\sico})$ is a separable
extension of~$F(G)$ of degree~$|\Pi_{G}|$.

\smallbreak

The same applies to the morphism
$\pi\times\id_{T}\colon G_{\sico}\times_{F}T\too G\times_{F}T$.
In particular, there exists a corestriction (or also called transfer) homomorphism
$$
\cor_{F(G_{\sico}\times_{F}T)/F(G\times_{F}T)}\, :\;\Brtor{n}(F(G_{\sico}\times_{F}T))
\,\too\,\Brtor{n}(F(G\times_{F}T))\, .
$$
By the projection formula the composition
$$
\cor_{F(G_{\sico}\times_{F}T)/F(G\times_{F}T)}\;\circ\;
\res_{F(G_{\sico}\times_{F}T)/F(G\times_{F}T)}
$$
is equal to the multiplication by
$[F(G_{\sico}\times_{F}T):F(G\times_{F}T)]=|\Pi_{G}|$
and so if the cardinality~$|\Pi_{G}|$ of the fundamental
group of~$G$ is coprime to~$n$ we get that the restriction
homomorphism~$\res_{F(G_{\sico}\times_{F}T)/F(G\times_{F}T)}$ is injective
on the subgroup of elements of exponent dividing~$n$.

\medbreak

We assume now that the integer~$n$ is coprime to~$|\Pi_{G}|$
and consider the commutative diagram
\begin{equation}
\label{GscDiag}
\xymatrix{
\Brtor{n}(F(G_{\sico}\times_{F}T)) & & \Brtor{n}(F(G\times_{F}T))
    \ar[ll]_-{\res_{F(G_{\sico}\times_{F}T)/F(G\times_{F}T)}}
\\
\Brtor{n}(G_{\sico}\times_{F}T) \ar[u]^-{\iota_{\sico}^{\ast}} & & 
      \Brtor{n}(G\times_{F}T) \ar[u]_-{\iota^{\ast}}
               \ar[ll]_-{(\pi\times\id_{T})^{\ast}}
\\
 & \Brtor{n}(T)\, , \ar[ru]_-{p_{T}^{\ast}} \ar[lu]^-{p_{\sico}^{\ast}}
}
\end{equation}
where~$\iota$ and~$\iota_{\sico}$ are the respective generic points,
and~$p_{\sico}$ and~$p_{T}$ denote the projection $G_{\sico}\times_{F}T\too T$
and $G\times_{F}T\too T$, respectively.

\smallbreak

Since both $G\times_{F}T$ and $G_{\sico}\times_{F}T$ are smooth
we know that~$\iota_{\sico}^{\ast}$~and~$\iota^{\ast}$ are
both injective by Auslander and Goldman~\cite[Thm.\ 7.2]{AuGo60},
and so by the commutative diagram above we conclude
that $(\pi\times\id_{T})^{\ast}:\Brtor{n}(G\times_{F}T)\too
\Brtor{n}(G_{\sico}\times_{F}T)$ is injective, too. Since as
shown above, see Corollary~\ref{cor34}, the
homomorphism~$p_{\sico}^{\ast}$ is an isomorphism this
implies by a diagram chase on~(\ref{GscDiag}) that
$p_{T}^{\ast}:\Brtor{n}(T)\too\Brtor{n}(G\times_{F}T)$ is surjective
and so also an isomorphism. We have proven:
\end{emptythm}

\begin{emptythm}
\label{SemisimpleGrThm}
{\bfseries Corollary.}
{\it
Let~$F$ be a field and~$n\geq 2$ an integer which
is not divisible by $\khar F$. Let further~$G$ be a
connected semisimple algebraic group
over a field~$F$ and~$T$ an $F$-scheme, such that~$T_{s}$
is isomorphic to a torus over~$F_{s}$. Assume that~$n$ is
coprime to the order of the fundamental group~$\Pi_{G}$
of~$G$. If moreover~$\khar F$ does not divide~$|\Pi_{G}|$
then the pull-back homomorphism
$$
p_{T}^{\ast}\, :\;\Brtor{n}(T)\,\too\,\Brtor{n}(G\times_{F}T)
$$
along the projection $p_{T}\colon G\times_{F}T\too T$
is an isomorphism.
}
\end{emptythm}

\begin{emptythm}
\label{GlnExplSubSect}
{\it An example.}
Let~$F$ be a field and~$\Gl_{d}(F)$ the general
linear group over~$F$ of rank~$d$. As an $F$-variety
(not as an algebraic group!) the general linear
group $\Gl_{d}(F)$ is isomorphic to~$\Sl_{d}(F)\times_{F}\Gl_{1}(F)$.
Since the special linear group is simply connected we get
from~\ref{scGrSubSect} an isomorphism
$$
\Brtor{n}(\Gl_{d}(F))\,\simeq\,\Brtor{n}(F[t,t^{-1}])
$$
for all~$n\in\N$ which are not divisible by~$\khar F$.
\end{emptythm}

\begin{emptythm}
\label{redGrSepclFieldSubSect}
{\it Brauer groups of reductive groups over a separably
closed field.}
Let~$F$ be a separably closed field whose characteristic
does not divide the integer~$n\geq 2$, and~$H$
a reductive group over~$F$. Then by~\cite[Lem.\ 8]{Ma78}
there exists a (split) torus $T\hookrightarrow H$, such
that we have an isomorphism
\begin{equation}
\label{algclosedDecEq}
H\,\simeq\, G\times_{F}T
\end{equation}
as $F$-schemes, where~$G$ is the commutator subgroup~$[H,H]$.
Let us recall the proof of
decomposition~(\ref{algclosedDecEq}). Let~$T_{1}$ be
a maximal torus of~$G=[H,H]$. It is contained in a
maximal torus~$S$ of~$H$, and since~$F$ is separably
closed we have $S=T_{1}\times_{F} T$ for some torus~$T\subset S$.
The centralizer of~$T_{1}$ in~$G$ is equal to~$T_{1}$ and
so we have $T(F_{s})\cap G(F_{s})=\{ e\}$, hence the
morphism of $F$-schemes $G\times_{F}T\too H$,
$(g,t)\mapsto g\cdot t$ is an isomorphism. (Note that
if~$H$ is semisimple, \ie $H=[H,H]$,
then the torus~$T$ in~(\ref{algclosedDecEq}) is
the ``trivial'' torus~$\Spec F$.)

\medbreak

If~$\khar F$ does not divide the order of the fundamental
group~$\Pi_{G}$ of~$G$ and the integer~$n$ is coprime
to this order, we conclude from this isomorphism
of $F$-schemes and Corollary~\ref{cor34} that
$\Brtor{n}(H)\simeq\Brtor{n}(T)$.
If the integer~$n$ is not coprime to~$|\Pi_{G}|$ we use
the strategy of Magid~\cite[Thm.\ 9]{Ma78} who has computed
the Brauer group of a reductive group over an algebraically
closed field of characteristic zero.
The same approach gives also some results for fields
of positive characteristic, and fields which are not
separably closed.
\end{emptythm}

\begin{emptythm}
\label{redGrArbitraryFieldSubSect}
{\it Some remarks about Brauer groups of reductive groups
over an arbitrary field.}
Let~$F$ be a field (not assumed to be separably closed)
and~$H$ a reductive group over~$F$. Let~$G=[H,H]$ and
assume the following:
\begin{itemize}
\item[(i)]
There is an isomorphism of $F$-schemes
$$
H\,\simeq\, G\times_{F}T\, ,
$$
where~$T$ is a split $F$-torus, and

\smallbreak

\item[(ii)]
the fundamental group~$\Pi_{G}$ of the semisimple
algebraic group $G=[H,H]$ is constant, \ie defined
over the base field, and $\khar F$ does not divide
the order of~$\Pi_{G}$. Hence the
simply connected cover $\pi:G_{\sico}\too G$
is a Galois cover with group~$G$.
\end{itemize}

\smallbreak

We have then a commutative diagram
$$
\xymatrix{
\Brtor{n}(G_{\sico}\times_{F}T) & & \Brtor{n}(G\times_{F}T)
  \ar[ll]_-{(\pi\times\id_{T})^{\ast}}
\\
 & \Brtor{n}(T) \rlap{\, } \ar[lu]^-{p_{\sico}^{\ast}} \ar[ru]_-{p_{T}^{\ast}}
}
$$
where $p_{\sico}:G_{\sico}\times_{F}T\too T$ and
$p_{T}:G\times_{F}T\too T$ are the respective projections.
Since as remarked above
$p_{\sico}^{\ast}:\Brtor{n}(T)\too \Brtor{n}(G_{\sico}\times_{F}T)$
is an isomorphism, and since~$p_{T}^{\ast}$ is a split
monomorphism, we get from this diagram an isomorphism
$$
\Brtor{n}(G\times_{F}T)\,\simeq\,\Brtor{n}(T)\,\oplus\,
\Ker\big(\,\Brtor{n}(G\times_{F}T)\,\xrightarrow{\; (\pi\times\id_{T})^{\ast}\;}\,
\Brtor{n}(G_{\sico}\times_{F}T)\,\big)\, .
$$

\smallbreak

By our assumptions
$\pi\times\id_{T}:G_{\sico}\times_{F}T\too G\times_{F}T$
is a Galois cover with group~$\Pi_{G}$. Hence,
see~\ref{HSSpSeqSubSect}, we have a convergent
spectral sequence
$$
E_{2}^{p,q}\, :=\;\HM^{p}(\Pi_{G},\HM^{q}_{\et}(G_{\sico}\times_{F}T,\Gm))\;\;
\Longrightarrow\;\HM^{p+q}_{\et}(G\times_{F}T,\Gm)\, ,
$$
where $\HM^{i}(\Pi_{G},\, -\, )$ denotes the group cohomology
of the finite group~$\Pi_{G}$. The associated five
term exacts sequence reads as follows:
$$
\Pic (G\times_{F}T)\,\too\,\Pic (G_{\sico}\times_{F}T)^{\Pi_{G}}\,\too\,
\HM^{2}(\Pi_{G},F[G_{\sico}\times_{F}T]^{\times})\,\too\,
\qquad\qquad\quad
$$
$$
\qquad\quad
\Ker\Big(\,\Br (G\times_{F}T)\,
\xrightarrow{\; (\pi\times\id_{T})^{\ast}\;}\,
\Br (G_{\sico}\times_{F}T)\,\Big)\,\too\,
\HM^{1}(\Pi_{G},\Pic (G_{\sico}\times_{F}T))\, .
$$
(This is the Chase-Harrison-Rosenberg~\cite{ChHaRo69} exact
sequence for the Galois extension $G_{\sico}\times_{F}T\too
G\times_{F}T$, see also~\cite[Chap.\ IV]{SAC}.)

\smallbreak

We have~$\Pic G_{\sico}=0$, see \eg~\cite[Thm.\ 1.5]{Gi09},
and so by localization and homotopy invariance also
$\Pic (G_{\sico}\times_{F}T)=0$ as~$T$ is a split torus
by our assumptions. Therefore we get
from this exact sequence an isomorphism
$$
\Brtor{n}(G\times_{F}T)\,\simeq\,\Brtor{n}(T)\,\oplus\,
\HMtor{n}^{2}(\Pi_{G},F[G_{\sico}\times_{F}T]^{\times})\, ,
$$
where $\HMtor{n}^{i}(\Pi_{G},\, -\, )$ denotes the
subgroup of $\HM^{i}(\Pi_{G},\, -\, )$ of elements
of order dividing~$n$.

\smallbreak

By Rosenlicht~\cite{Ros61}
we have $F[G_{\sico}\times_{F}T]^{\times}= F^{\times}\times\widehat{T}$,
where $\widehat{T}=\Hom (T,\Gm)$ is the group of characters of~$T$.
The group~$\Pi_{G}$ operates trivially on~$\widehat{T}$, and therefore
there is an isomorphism
$$
\HMtor{n}^{2}(\Pi_{G},F[G_{\sico}\times_{F}T]^{\times})\,\simeq\,
\HMtor{n}^{2}(\Pi_{G},F^{\times})\,\oplus\,\HMtor{n}^{2}(\Pi_{G},\Z)^{\oplus\, d}\, ,
$$
where~$d$ is the rank of the torus~$T$. Since~$\Pi_{G}$ is
a finite group operating trivially on~$\Z$ we have
$\HM^{1}(\Pi_{G},\Z)=0$ and so by the long exact cohomology
sequence for the short exact sequence
$\Z\xrightarrow{\cdot n}\Z\too\Z /\Z n$ we get
$\HMtor{n}^{2}(\Pi_{G},\Z)\simeq\HM^{1}(\Pi_{G},\Z/\Z n)$.

\smallbreak

We arrive finally at the following result, which is due
to Magid~\cite[Thm.\ 9]{Ma78} if the field~$F$ is algebraically closed of
characteristic zero.
\end{emptythm}

\begin{emptythm}
\label{sepClosedThm}
{\bfseries Theorem.}
{\it
Let $F$ be a  field and~$H$ a (connected) reductive
$F$-group. Assume that~$\khar F$ neither divides~$n$ nor the
order of the the fundamental group~$\Pi_{G}$ of the
semisimple $F$-group $G=[H,H]$. Assume further that
\begin{itemize}
\item[(i)]
there is an isomorphism of $F$-schemes
$H\simeq G\times_{F}T$, where~$T$ is a
split $F$-torus of rank $d$, and

\smallbreak

\item[(ii)]
the simply connected cover $\pi:G_{\sico}\too G$
is a Galois cover with group~$\Pi_{G}$.
\end{itemize}

\smallbreak

Then we have an isomorphism
$$
\Brtor{n}(H)\,\simeq\,\Brtor{n}(T)\,\oplus\,\HMtor{n}^{2}(\Pi_{G},F^{\times})
\,\oplus\,\HM^{1}(\Pi_{G},\Z/\Z n)^{\oplus\, d}\, .
$$
}
\end{emptythm}
\medbreak

\begin{emptythm}
{\bfseries Remark.}
Assumptions~(i) and~(ii) are automatically satisfied if
the field~$F$ is separably closed.
\end{emptythm}

\begin{emptythm}
\label{redGrAdjSubSect}
Let again~$H$ be a connected reductive linear algebraic
group over a field~$F$, which is not assumed to be
separably closed. Let further
(as above)~$n$ be an integer which is coprime to the
characteristic of~$F$. We denote by~$\bar{F}$ an algebraic
closure of~$F$.

\smallbreak

By~\cite[\S 2]{BoTi65} the derived group as well as
the radical of~$H$ are defined over~$F$ and we have
$$
H(\bar{F})\, =\, [H,H](\bar{F})\cdot\Rad (H)(\bar{F})\, ,
$$
where~$[H,H]$ and~$\Rad (H)$ denote the $F$-subgroup
schemes of~$H$ which represent the commutator subgroup and
the radical of~$H$, respectively.

\smallbreak

Moreover, the radical is a torus and the intersection
of the $\bar{F}$-points of~$[H,H]$ and~$\Rad (H)$ is
a finite group. Since this intersection is in the center
of~$H$ the group~$H$ is isomorphic to $[H,H]\times_{F}\Rad (H)$
if~$[H,H]$ is of adjoint type. We conclude from
Corollary~\ref{SemisimpleGrThm}:

\medbreak

\begin{emptythm}
\label{cor312}
{\bfseries Corollary.}
{\it
Let~$H$ be a connected reductive group over a field~$F$,
such that~$[H,H]$ is of adjoint type. Assume that
the characteristic of~$F$ does not divide the order
of the fundamental group of the commutator subgroup of~$H$.
Then
$$
\Brtor{n}(\Rad (H))\,\simeq\,\Brtor{n}(H)
$$
for all integers~$n$ which are not divisible
by~$\khar F$ and coprime to the order of the
fundamental group of the commutator subgroup~$[H,H]$.
}
\end{emptythm}
\end{emptythm}

\goodbreak
\section{The Brauer group of geometrically integral and regular
rings of finite type over a field of positive characteristic}
\label{posCharSect}\bigbreak

\begin{emptythm}
\label{polynomialringThm}
{\it The Brauer group of polynomial rings over fields of positive
characteristic.}
We start with the following result due to Auslander and Goldman~\cite{AuGo60}
if~$d=1$, and due to Knus, Ojanguren and Saltman~\cite{KnOjSa76} if~$d\geq 2$.

\medbreak

\begin{emptythm}
{\bfseries Proposition.}
{\it
Let~$F$ be a field of characteristic~$p>0$. Then if~$d\geq 2$,
or if~$d\geq 1$ and~$F$ is not perfect, the natural split monomorphism
$$
\Br (F)\,\too\,\Br (F[T_{1},\ldots ,T_{d}])
$$
is not surjective.
}
\end{emptythm}
\begin{proof}
If~$F$ is not perfect we know by
Knus, Ojanguren and Saltman~\cite[Thm.\ 5.5]{KnOjSa76}
that~$\Br (F[T_{1}])\simeq\Br (F)\oplus\Z (p^{\infty})^{(I)}$
for some infinite set~$I$, see also Auslander and
Goldman~\cite[Thm.\ 7.5]{AuGo60}. Since~$\Br (F[T_{1}])$ is a direct
summand of~$\Br (F[T_{1},\ldots ,T_{d}])$ the claim follows.


Assume now that~$F$ is perfect and~$d\geq 2$.
If~$F$ is moreover finite then this result has
been proven by Knus, Ojanguren, and Saltman
in~\cite[Proof of Thm.\ 5.7]{KnOjSa76}. The same
argument works for any perfect field. In fact,
since~$F$ is perfect we have
$\big( F[T_{1}^{\frac{1}{p}},\ldots ,T_{d}^{\frac{1}{p}}]\big)^{p}=
F[T_{1},\ldots ,T_{d}]$, and so by~\cite[Cor.\ 3.10]{KnOjSa76}
the subgroup $\Brtor{p}(F[T_{1},\ldots ,T_{d}])$ is isomorphic
to the kernel of the natural homomorphism
$\Br (F[T_{1},\ldots ,T_{d}])\too
\Br (F[T_{1}^{\frac{1}{p}},\ldots ,T_{d}^{\frac{1}{p}}])$. This
kernel contains the kernel of $\Br (F[T_{1},\ldots ,T_{d}])\too
\Br (F[T_{1},\ldots ,T_{n-1},T_{d}^{\frac{1}{p}}])$ which in turn
contains a subgroup isomorphic to the quotient group
$F[T_{1},\ldots ,T_{d-1}]\big/(F[T_{1},\ldots ,T_{d-1}])^{p}$
by~\cite[Prop.\ 5.3]{KnOjSa76}. The latter group is not
trivial since by assumption~$n\geq 2$, and so also
$\Brtor{p}(F[T_{1},\ldots ,T_{d}])\not= 0$.

On the other hand, by a result of Albert~\cite[p.\ 109]{StAlg}
any central simple $F$-algebra of exponent~$p$ is Brauer equivalent
to a cyclic algebra and so trivial in~$\Br (F)$ since~$F$ is perfect.
Therefore we have $\Brtor{p}(F)=0$. We are done.
\end{proof}
\end{emptythm}

\begin{emptythm}
\label{connectedSubSect}
Our aim is now to show that if~$F$ is not perfect, or if~$n\geq 2$,
then the cokernel of the pull-back homomorphism
$\Br (F)\too\Br (F[T_{1},\ldots ,T_{d}])$
is not empty, $p$-torsion, and $p$-divisible for all
fields~$F$ of characteristic~$p>0$. We prove first
the following general lemma;
see also Orzech and Small~\cite[Proof of Cor.\ 8.8]{BrCR}.

\medbreak

\begin{emptythm}
{\bfseries Lemma.}
{\it
Let~$R$ be a regular integral domain of characteristic~$p>0$
and~$l>0$ an integer which is not divisible by~$p$. Then the
pull-back
$$
\pi_{R}^{\ast}\, :\;\Brtor{l} (R)\,\too\,\Brtor{l} (R[T])
$$
along the inclusion $\pi_{R}:R\too R[T]$
is an isomorphism. In particular, if~$F$ is a field of
characteristic~$p$ and~$p$ does not divide~$l$ then the
pull-back
$$
\Brtor{l}(F)\too\Brtor{l}(F[T_{1},\ldots ,T_{d}])
$$
is an isomorphism.
}
\end{emptythm}
\begin{proof}
We denote for a ring~$S$ by~$a_{S}$ the zero section
$S[T]\too S$, $T\mapsto 0$. This is a left inverse of
the inclusion $\pi_{S}:S\too S[T]$, \ie~$a_{S}\circ\pi_{S}=\id_{S}$.

We have a commutative diagram
$$
\xymatrix{
\Brtor{l}(R[T]) \ar[r]^-{a_{R}^{\ast}} \ar[d] & \Brtor{l}(R) \ar[d]
\\
\Brtor{l} (K[T]) \ar[r]^-{a_{K}^{\ast}} & \Brtor{l}(K) \rlap{\, ,}
}
$$
where~$K$ is the quotient field of~$R$ and the vertical
arrows are induced by the inclusions~$R\hookrightarrow K$ and
$R[T]\hookrightarrow K[T]$.

\smallbreak

By a result of Grothendieck, see \eg~\cite[Lem.\ 4.4]{Gi09},
the pull-back homomorphism~$\pi_{K}^{\ast}:\Brtor{l}(K)\too\Brtor{l}(K[T])$
is an isomorphism since (by assumption)~$p=\khar K$ does
not divide~$l$. Hence
$a_{K}^{\ast}= (\pi_{K}^{\ast})^{-1}:\Brtor{l}(K[T])\too\Brtor{l}(K)$
is also an isomorphism. Since~$R$ is a regular integral domain
we know from Auslander and Goldman~\cite[Thm.\ 7.2]{AuGo60}
that~$\Br (R)\too\Br (K)$ is injective and so above
commutative diagram shows that
$a_{R}^{\ast}:\Brtor{l}(R[T])\too\Brtor{l}(R)$ is a monomorphism.
This implies the lemma since $a_{R}^{\ast}\circ\pi_{R}^{\ast}$ is
equal to the identity of~$\Brtor{l}(R)$.
\end{proof}

We are in position to prove the following (slight)
generalization of Knus, Ojanguren, and
Saltman~\cite[Thm.\ 5.7]{KnOjSa76}.
\end{emptythm}

\begin{emptythm}
\label{BrauerGrPolynomialRingThm}
{\bfseries Theorem.}
{\it
Let~$F$ be a field of characteristic~$p>0$ and~$d\geq 1$ an integer.
Then we have
$$
\Br (F[T_{1},\ldots ,T_{d}])\,\simeq\,\Br (F)\,\oplus\,\Z (p^{\infty})^{(I)}\, ,
$$
where the set~$I$ is non empty if (and only if)~$F$
is not perfect or~$d\geq 2$.
}

\begin{proof}
We show first that
$\Coker\Big(\Br (F)\,\too\,\Br (F[T_{1},\ldots ,T_{d}])\Big)$
is $p$-torsion. Let for this~$0\not= x\in\Br (F[T_{1},\ldots ,T_{d}])$.
Since the Brauer group of a commutative ring is torsion, see
\eg~\cite[Chap.\ 12]{BrCR}, there exists~$d\geq 2$,
such that $d\cdot x=0$. Write $d=p^{m}\cdot d'$ with~$d'$
an integer coprime to~$p$ and~$m\geq 0$. Then $d'\cdot (p^{m}\cdot x)=0$
and so~$p^{m}\cdot x$ has exponent~$d'$. Since~$p$ does not
divide~$d'$ it follows from the lemma above
that~$p^{m}\cdot x$ is in the image of natural
homomorphism~$\Br (F)\too\Br (F[T_{1},\ldots ,T_{d}])$.

\smallbreak

But~$\Br (F[T_{1},\ldots ,T_{d}])$ is $p$-divisible
by~\cite[Cor.\ 4.4]{KnOjSa76} and so is also the
quotient $\Br (F[T_{1},\ldots ,T_{d}])\big/\Br (F)$.
Therefore by the structure theorem for divisible
abelian groups, see~\cite[Thm.\ 4]{IAb},
it is isomorphic to~$\Z (p^{\infty})^{(I)}$ for some set~$I$. By
Theorem~\ref{polynomialringThm} this set is not empty if~$d\geq 1$
and~$F$ is not perfect, or if~$d\geq 2$. We are done.
\end{proof}
\end{emptythm}

\begin{emptythm}
\label{NormalizationSubSect}
Let~$F$ be a field of characteristic~$p>0$
and~$R$ a regular $F$-algebra of finite type. We assume that~$R$
is geometrically integral and has Krull dimension~$d\geq 1$.
Then, see \eg~\cite[Cor.\ 16.18]{ComAlg},
there are~$T_{1},\ldots ,T_{d}\in R$, such that
\begin{itemize}
\item[(i)]
$R$ is finite over the (polynomial) ring~$D:=F[T_{1},\ldots ,T_{d}]$, and

\smallbreak

\item[(ii)]
the quotient field~$K$ of~$R$ is a separable extension
of~$E:=F(T_{1},\ldots ,T_{d})$.
\end{itemize}
(In other words, $T_{1},\ldots ,T_{d}\in R$ are a separating transcendence basis
of~$K$ over~$F$.)

\smallbreak

We set now $\bBr (S):=\Coker\Big(\,\Br (F)\,\too\,\Br (S)\,\Big)$
for an $F$-algebra~$S$.
With this notation we have a commutative diagram where all
morphisms are induced by the ring inclusions except for~$c_{K/E}$
which is induced by the corestriction
\begin{equation}
\label{NormalizationEq}
\xymatrix{
\bBr (R) \ar[r]^-{\iota_{R}} & \bBr (K) \ar@/^/[d]^-{c_{K/E}}
\\
\bBr (D) \ar[u]^-{r_{R/D}} \ar[r]_-{\iota_{D}} & \bBr (E) \ar@/^/[u]^-{r_{K/E}}
\rlap{\, .}
}
\end{equation}
Since~$R$ and~$D$ are regular integral domains the pull-back
maps $\Br (R)\too\Br (K)$ and $\Br (D)\too\Br (E)$
are monomorphisms by Auslander and Goldman~\cite[Thm.\ 7.2]{AuGo60}
and so also~$\iota_{R}$ and~$\iota_{D}$ are both injective.

\medbreak

Assume now that~$\dim R\geq 1$ and~$F$ is not perfect,
or that~$\dim R\geq 2$. Then by Theorem~\ref{BrauerGrPolynomialRingThm}
we have
$$
\bBr (D)\, =\,\bBr (F[T_{1},\ldots ,T_{n}])\,\simeq\,\Z (p^{\infty})^{(I)}
$$
for some non empty set~$I$. Since~$K\supseteq E$ is a finite separable
extension we know that $c_{K/E}\circ r_{K/E}$ is equal to the multiplication
by~$m=[K:E]$. Therefore by the commutative diagram~(\ref{NormalizationEq})
we have
$$
(c_{K/E}\circ\iota_{R})\big[\,r_{R/D}(\bBr (D))\,\big]\, =\,
m\cdot \iota_{D}(\bBr (D))\,\simeq\,\Z (p^{\infty})^{(I)}
$$
for some non empty set~$I$ since~$\Z (p^{\infty})$
is divisible. Hence~$\bBr (R)$ is not trivial.

\smallbreak

We have shown the following theorem.
\end{emptythm}

\begin{emptythm}
\label{mainp-Thm}
{\bfseries Theorem.}
{\it
Let~$F$ be a field of characteristic~$p>0$.
If~$R$ is a regular and geometrically integral $F$-algebra
of dimension~$\geq 1$ if~$F$ is non perfect, respectively
of dimension~$\geq 2$ if~$F$ is perfect, then the natural
homomorphism
$$
\Br (F)\,\too\Br (R)
$$
is not surjective.
}
\end{emptythm}

\begin{emptythm}
\label{scExpl}
{\bfseries Examples.}
Let~$F$ be a field of characteristic~$p>0$.
\begin{itemize}
\item[(i)]
Let~$G$ be a (connected) simply connected
semisimple linear algebraic group over~$F$.
Since~$G(F)\not=\emptyset$ we have
$\Br (G)\,\simeq\,\Br (F)\,\oplus\, M$.
for some subgroup~$M$ which is not trivial
by Theorem~\ref{mainp-Thm} above and $p$-divisible
by~\cite[Cor.\ 4.4]{KnOjSa76}.
In~\cite{Gi09} it has been shown that the pull-back
homomorphism $\Brtor{n}(F)\too\Brtor{n}(G)$ is an
isomorphism for all~$n\in\N$ which are not
divisible by $p=\khar F$ and so, \cf the proof
of Theorem~\ref{BrauerGrPolynomialRingThm}, the
subgroup~$M$ is also $p$-torsion. Hence we have
$$
\Br (G)\,\simeq\,\Br (F)\,\oplus\,\Z (p^{\infty})^{(I)}
$$
for some non empty set~$I$.

\medbreak

\item[(ii)]
We assume here that $p\not= 2$. Let~$X_{q,\aff}$ be a non
singular affine quadric over~$F$ as in~\ref{affQuadSubSect}.
We assume that~$X_{q,\aff}(F)\not=\emptyset$. Then
by Theorem~\ref{affQuadThm} the canonical homomorphism
$\Brtor{n}(F)\too\Brtor{n}(X_{q,\aff})$ is an isomorphism
for all~$n\in\N$ not divisible by~$\khar F$. Therefore
by the same arguments as in part~(i) above we have
$\Br (X_{q,\aff})\simeq\Br (F)\oplus\Z (p^{\infty})^{(I)}$
for some non empty set~$I$.
\end{itemize}
\end{emptythm}

\begin{emptythm}
\label{globalRem}
{\bfseries Remark.}
Theorem~\ref{mainp-Thm} is wrong for non affine
schemes. For instance, a projective homogeneous
variety~$X$ for a semisimple algebraic group over
a field~$F$ (\eg a projective quadric) is smooth
and geometrically integral, but as shown by Merkurjev
and Tignol~\cite[Thm.\ B]{MeTi95} the natural map
$\Br (F)\too\HM^{2}_{\et}(X,\Gm)$ is surjective even
if $\khar F>0$.
\end{emptythm}

\bibliographystyle{amsalpha}

\end{document}